\documentclass[a4paper,11pt,reqno]{amsart}

\usepackage{dsfont,amssymb}

\usepackage[all]{xy}

\numberwithin{equation}{section}

\newtheorem{theorem}{Theorem}[section]
\newtheorem{corollary}[theorem]{Corollary}
\newtheorem{lemma}[theorem]{Lemma}
\newtheorem{proposition}[theorem]{Proposition}

\theoremstyle{remark}
\newtheorem{remark}[theorem]{Remark}
\newtheorem{example}[theorem]{Example}

\theoremstyle{definition}
\newtheorem{definition}[theorem]{Definition}

\newcommand\bp{\begin{proof}}
\newcommand\ep{\end{proof}}
\newcommand\ee{\nopagebreak\mbox{\ }\hfill$\diamondsuit$}

\newcommand\un{{\mathds 1}}

\newcommand\Dhat{{\hat\Delta}}

\newcommand\End{\operatorname{End}}

\newcommand\Hom{\operatorname{Hom}}
\newcommand\Hilb{\operatorname{Corr}}
\newcommand\Hilbf{\operatorname{Hilb}_f}

\newcommand\Mor{\operatorname{Mor}}
\newcommand\Rep{\operatorname{Rep}}
\newcommand\Tr{\operatorname{Tr}}

\newcommand{\C}{{\mathbb C}}

\newcommand\T{{\mathbb T}}
\newcommand\Z{{\mathbb Z}}

\newcommand{\CC}{{\mathcal C}}
\newcommand{\B}{{\mathcal B}}
\newcommand{\D}{{\mathcal D}}

\newcommand{\M}{{\mathcal M}}

\newcommand\U{{\mathcal U}}

\newcommand\eps{\varepsilon}

\begin{document}

\title{Duality theory for nonergodic actions}

\author[S. Neshveyev]{Sergey Neshveyev}

\email{sergeyn@math.uio.no}

\address{Department of Mathematics, University of Oslo,
P.O. Box 1053 Blindern, NO-0316 Oslo, Norway}

\thanks{The research leading to these results has received funding from the European Research Council
under the European Union's Seventh Framework Programme (FP/2007-2013) / ERC Grant Agreement no. 307663
}

\date{March 25, 2013; minor changes April 3, 2013}

\begin{abstract}
Generalizing work by Pinzari and Roberts, we characterize actions of a compact quantum group $G$ on C$^*$-algebras in terms of what we call weak unitary tensor functors from $\Rep G$ into categories of C$^*$-correspondences. We discuss the relation of our construction of a C$^*$-algebra from a functor to some well-known crossed product type constructions, such as cross-sectional algebras of Fell bundles and crossed products by Hilbert bimodules. We also relate our setting to recent work of De Commer and Yamashita by showing that any object in a module C$^*$-category over $\Rep G$ produces a weak unitary tensor functor, and, as a consequence, actions can also be described in terms of $(\Rep G)$-module C$^*$-categories. As an application we discuss deformations of C$^*$-algebras by cocycles on discrete quantum groups.
\end{abstract}

\maketitle

\section*{Introduction}

Category theory has, from the early beginning, played an important role in quantum groups. In the operator algebraic approach to quantum groups the key result connecting the two areas is due to Woronowicz~\cite{WorTK}. Generalizing the classical Tannaka-Krein duality he showed that by associating to a compact quantum group its representation category together with the canonical fiber functor, we get a duality between compact quantum groups on one side, and C$^*$-tensor categories with conjugates and a unitary fiber functor, on the other. Therefore, in principle, all properties of a compact quantum group $G$ can be formulated entirely in terms of its representation category $\Rep G$ and canonical fiber functor. Remarkably, a lot of properties depend on $\Rep G$ alone. A systematic study of such properties was made possible by Bichon, De Rijdt and Vaes~~\cite{BDRV}, who showed how, given two monoidally equivalent (that is, having equivalent representation categories) compact quantum groups, to construct a linking C$^*$-algebra connecting the two. The linking algebra is equipped with ergodic actions of both quantum groups, and by fixing one quantum group $G$ and varying the other, or in other words, by considering all possible unitary fiber functors on $\Rep G$, we get all ergodic actions of $G$ of full quantum multiplicity. This categorical point of view on actions, together with the construction of linking algebra, has been extremely successful. It has been applied to a variety of seemingly unrelated problems, from a study of random walks on discrete quantum groups~\cite{DRV} to $K$-theoretic computations~\cite{Voigt}.

The results of Bichon, De Rijdt and Vaes were later generalized by Pinzari and Ro\-berts~\cite{PR}, who showed how to describe all ergodic actions of $G$ in terms of $\Rep G$. Namely, for every ergodic action they constructed a ``spectral functor'' from $\Rep G$ into the category of Hilbert spaces, and then gave an abstract characterization of such functors. Their result implies, in particular, that isomorphism classes of ergodic actions of monoidally equivalent compact quantum groups are in canonical correspondence with each other. Soon afterwards De Rijdt and Vander Vennet~\cite{DRV} showed that the same is true even for nonergodic actions. Their argument bypasses category theory altogether and is based on induction using the linking algebra. A natural problem completing this circle of ideas is nevertheless to find a description of actions of $G$ entirely in terms of~$\Rep G$. Our goal in this paper is to do exactly that. By modifying the definition of a spectral functor and the axioms of Pinzari and Roberts, we show that actions of a compact quantum group~$G$ correspond to a class of functors, which we call weak unitary tensor functors, from $\Rep G$ into categories $\Hilb A$ of C$^*$-correspondences over C$^*$-algebras $A$. It should become apparent from our results, and is not difficult to show directly, that in the case $A=\C$ our definitions/axioms are equivalent to the ones given by Pinzari and Roberts. Overall, the construction of a C$^*$-algebra from a functor $\Rep G\to\Hilb A$ follows familiar lines going back to Woronowicz~\cite{WorTK}. Since some of the maps involved are not adjointable, we just have to be more careful not to overuse various dualities.

A different solution to the same problem is suggested by recent work of De Commer and Yamashita~\cite{DCY}. Complementing the results of Pinzari and Roberts, they showed that ergodic actions of~$G$ can also be described in terms of semisimple $(\Rep G)$-module C$^*$-categories with a fixed simple generating object. In fact, a significant part of their arguments does not involve ergodicity/semisimplicity in any way, and we show that indeed by discarding these assumptions we get a characterization of general actions in terms of module categories. The relation between two categorical pictures can be described as follows. Given a right $(\Rep G)$-module C$^*$-category and an object $M$ in it, the functor $\Mor(M,M\otimes\cdot)$ has a canonical structure of a weak unitary tensor functor. Therefore, using the analogy with representation theory, we can say that the relation between two pictures is as between a cyclic representation and its matrix coefficient defined by the cyclic vector. From this point of view weak unitary tensor functors are categorifications of positive definite functions. It is interesting that the module category approach, being less economical than the approach via weak tensor functors, seems, nevertheless, more suitable for classification of actions, at least for representation categories described by simple universal properties~\cite{DCY2}.

In the last section of this paper we discuss some examples and applications of our general results. The construction of a C$^*$-algebra from a weak unitary tensor functor is reminiscent of various crossed product type constructions. To make the connection more explicit, we reformulate this construction in a category-free way. This will make it clear that for duals of discrete groups it generalizes such constructions as cross-sectional algebras of Fell bundles or crossed products by Hilbert bimodules. We also show that categorical point of view on actions naturally leads to a construction of deformation of C$^*$-algebras by  $2$-cocycles on discrete quantum groups.

\bigskip

\section{Spectral functors} \label{s1}

Let us first fix the notation; overall, we follow the same conventions as in \cite{NTbook}. Consider a compact quantum group $G$. The Hopf $*$-algebra of matrix coefficients of finite dimensional representations of $G$ is denoted by $(\C[G],\Delta)$. A finite dimensional representation $U$ of $G$ is an invertible element of $B(H_U)\otimes C(G)$ such that $(\iota\otimes\Delta)(U)=U_{12}U_{13}$. The tensor product of two representations $U$ and $V$ is denoted by $U\times V$ and is defined by $U\times V=U_{13}V_{23}$.

The contragredient representation to a finite dimensional representation $U$ is defined by
$$
U^c=(j\otimes\iota)(U^{-1})\in B(H_U^*)\otimes C(G),
$$
where $j$ is the canonical anti-isomorphism $B(H_U)\cong B(H_U^*)$. When $H_U$ is a Hilbert space, we identify the dual space $H_U^*$ with the complex conjugate Hilbert space $\bar H_U$.

We denote the Woronowicz character $f_1\in\C[G]^*$ by $\rho$. For every finite dimensional representation~$U$ of $G$ we have a representation $\pi_U$ of the algebra $\C[G]^*$ on $H_U$ defined by $\pi_U(\omega)=(\iota\otimes\omega)(U)$. Given a finite dimensional unitary representation $U$ of $G$, the conjugate representation is defined by
$$
\bar U=(j(\pi_U(\rho)^{1/2})\otimes1)U^c(j(\pi_U(\rho)^{-1/2})\otimes1)\in B(\bar H_U)\otimes C(G).
$$
This is a unitary representation equivalent to $U^c$, and $\pi_{\bar U}(\rho)=j(\pi_U(\rho)^{-1})$. Unitarity of $\bar U$ essentially characterizes $\rho$: if $U$ is an irreducible unitary representation, then $\pi_U(\rho)$ is the unique strictly positive operator in $B(H_U)$ such that the above definition of $\bar U$ gives a unitary element and such that $\Tr(\pi_U(\rho))=\Tr(\pi_U(\rho)^{-1})$. We will usually suppress $\pi_U$ and simply write $\rho\xi$ for $\xi\in H_U$ instead of~$\pi_U(\rho)\xi$.

\medskip

Denote by $\Rep G$ the C$^*$-tensor category of finite dimensional unitary representations of $G$. In this category $\bar U$ is conjugate to $U$, in the sense that there exist morphisms $R_U\colon\un\to \bar U\times U$ and $\bar R_U\colon\un\to U\times \bar U$, where $\un$ is the trivial representation of $G$ on the one-dimensional space $\C$, such that the compositions
$$
U \xrightarrow{\iota \otimes R_U} U\otimes \bar U \otimes U \xrightarrow{\bar R^*_U \otimes \iota} U \quad \mbox{and} \quad \bar U \xrightarrow{\iota \otimes \bar R_U} \bar U \otimes U \otimes \bar U \xrightarrow{R^*_U \otimes \iota}\bar U
$$
are the identity morphisms. Using the Woronowicz character $\rho$ we can define such morphisms by
\begin{equation}\label{econj}
R_U(1)=\sum_i\bar\xi_i\otimes\rho^{-1/2}\xi_i,\ \ \bar R_U(1)=\sum_i\rho^{1/2}\xi_i\otimes\bar\xi_i,
\end{equation}
where $\{\xi_i\}_i$ is an orthonormal basis in $H_U$. Note that the above expressions do not depend on the choice of an orthonormal basis, and
$$
\|R_U\|=\|\bar R_U\|=(\dim_q U)^{1/2},
$$
where $\dim_qU=\Tr\pi_U(\rho)$ is the quantum dimension of $U$.

\medskip

Consider now a continuous left action $\theta$ of $G$ on a C$^*$-algebra $B$, so $\theta\colon B\to C(G)\otimes B$ is an injective $*$-homomorphism such that $(\Delta\otimes\iota)\theta=(\iota\otimes\theta)\theta$ and $(C(G)\otimes1)\theta(B)$ is dense in $C(G)\otimes B$. Consider the $*$-subalgebra $\B\subset B$ consisting of elements $x\in B$ such that $\theta(x)$ lies in the algebraic tensor product $\C[G]\otimes B$. Equivalently, $\B$ is the linear span of elements of the form $(h\otimes\iota)((a\otimes1)\theta(x))$, where $a\in \C[G]$, $x\in B$ and $h$ is the Haar state on $G$. We call $\B$ the algebra of regular elements in~$B$. It is a dense $*$-subalgebra of $B$, and $\theta$ defines a left  coaction of the Hopf $*$-algebra $(\C[G],\Delta)$ on it. As follows from \cite[Lemma~2.5]{DCY0}, the positive map $E=(h\otimes\iota)\theta\colon B\to B^G$ is faithful on $\B$, in the sense that $E(x^*x)\ne0$ for every nonzero $x\in\B$.

Conversely, assume we have a left  coaction $\theta$ of the Hopf $*$-algebra $(\C[G],\Delta)$ on a $*$-algebra $\B$. By slightly extending the definition in \cite{DCY} we say that $\theta$ is an algebraic action of $G$ if the following conditions are satisfied:
\begin{itemize}
\item[(i)] the fixed point algebra $A=\B^G=\{x\in\B\mid\theta(x)=1\otimes x\}$ is a C$^*$-algebra;
\item[(ii)] the projection $E=(h\otimes\iota)\theta\colon\B\to A$ is positive and faithful, so $E(x^*x)\ge0$ and $E(x^*x)\ne0$ for $x\ne0$;
\item[(iii)] $E(x^*a^*ax)\le\|a\|^2E(x^*x)$ for all $a\in A$ and $x\in\B$.
\end{itemize}
Note that condition (iii) follows from (i) and (ii) if $\B$ is unital with unit $1\in A$. Note also that conditions (ii) and (iii) can be formulated by saying that $\B$ is a right pre-Hilbert $A$-module with inner product $\langle x,y\rangle=E(x^*y)$, and the operators of multiplication on the left by elements of $A$ are bounded.

Under the above conditions (i)-(iii) it is not difficult to show that the $*$-algebra $\B$ admits a unique C$^*$-completion $B$ such that $\theta$ extends to a continuous left action of the reduced form of $G$ on $B$, see \cite[Proposition~4.4]{DCY}. Namely, $\B$ is faithfully represented by operators of multiplication on the left on the right pre-Hilbert $A$-module~$\B$ with inner product $\langle x,y\rangle=E(x^*y)$, and this defines a norm on~$\B$. Note that in general the subalgebra of regular elements in the completion $B$ of $\B$ can be strictly larger than~$\B$.

\medskip

Given a finite dimensional unitary representation $U$ of $G$, we can consider $H_U$ as a left comodule over $(\C[G],\Delta)$ by defining
$$
\delta_U\colon H_U\to\C[G]\otimes H_U\ \ \text{by}\ \ \delta_U(\xi)=U^*_{21}(1\otimes\xi).
$$
Then, if $\theta$ is a continuous left action of $G$ on a C$^*$-algebra $B$, we can consider comodule maps $H_U\to B$. The linear span of images of all such maps is denoted by $B(U)$ and is called the spectral subspace of $B$ corresponding to $U$. Choosing representatives $U_\alpha$ of isomorphism classes of irreducible unitary representations of $G$, for the subalgebra $\B\subset B$ of regular elements we get
$$
\B=\bigoplus_\alpha B(U_\alpha).
$$
Consider the tensor product comodule $H_U\otimes\B$.
We denote by $(H_U\otimes B)^G$ the subcomodule of invariant vectors, so
$$
(H_U\otimes B)^G=\{X\in H_U\otimes B\mid U_{12}X_{13}=(\iota\otimes\theta)(X)\}.
$$
In other words, if $\{\xi_i\}_i$ is an orthonormal basis in $H_U$ and $U=(u_{ij})_{i,j}$ is written in the matrix form with respect to this basis, then
$$
(H_U\otimes B)^G=\{X=\sum_i\xi_i\otimes x_i\mid\theta(x_i)=\sum_ju_{ij}\otimes x_j\ \ \text{for all}\ \ i\}.
$$
Note that using Frobenius reciprocity we can identify $(H_U\otimes B)^G$ with $\Hom_G(H_{\bar U},\B)$, but we are not going to do this. The spectral subspaces can be recovered from $(H_U\otimes B)^G$ using the canonical surjective maps
$$
\bar H_U\otimes (H_U\otimes B)^G\to B(\bar U),\ \ \bar\xi\otimes X\mapsto(\bar\xi\otimes\iota)(X),
$$
which are isomorphisms for irreducible $U$.

The spaces $(H_U\otimes B)^G$ is our main object of interest. Clearly, if $A=B^G$, then these spaces are $A$-bimodules. Furthermore, if $X=\sum_i\xi_i\otimes x_i$ and $Y=\sum_i\xi_i\otimes y_i$ are vectors in $(H_U\otimes B)^G$, then the element $\sum_i x_i^*y_i$ is $G$-invariant. Hence $(H_U\otimes B)^G$ is a right Hilbert $A$-module with inner product $\langle X,Y\rangle=\sum_i x_i^*y_i$. This inner product is independent of the choice of an orthonormal basis, and by slightly abusing notation it can be written as $\langle X,Y\rangle=X^*Y$.

Given two finite dimensional unitary representation $U$ and $V$ of $G$, we have a map
$$
(H_U\otimes B)^G\otimes (H_V\otimes B)^G\to (H_{U\times V}\otimes B)^G, \ \ X\otimes Y\mapsto X_{13}Y_{23}.
$$
In other words, if we fix orthonormal bases $\{\xi_i\}_i$ in $H_U$ and $\{\zeta_j\}_j$ in $H_V$, then for $X=\sum_i\xi_i\otimes x_i$ and $Y=\sum_j\zeta_j\otimes y_j$ we have
$$
X\otimes Y\mapsto\sum_{i,j}\xi_j\otimes\zeta_j\otimes x_iy_j.
$$
It is obvious that this map defines an isometric map $(H_U\otimes B)^G\otimes_A (H_V\otimes B)^G\to (H_{U\times V}\otimes B)^G$.

\medskip

We are now ready to define spectral functors.

\begin{definition}
Given a continuous left action of a compact quantum group $G$ on a C$^*$-algebra $B$ with fixed point algebra $A$, the associated spectral functor is the unitary functor $F$ from $\Rep G$ into the C$^*$-tensor category $\Hilb A$ of C$^*$-correspondences over $A$ defined by
$$
F(U)=(H_U\otimes B)^G\ \ \text{with inner product}\ \ \langle X,Y\rangle=X^*Y
$$
for representations $U$, and $F(T)=T\otimes\iota$ for morphisms, together with the $A$-bilinear isometries
$$
F_{2,U,V}\colon F(U)\otimes_A F(V)\to F(U\times V),\ \ X\otimes Y\mapsto X_{13}Y_{23}.
$$
\end{definition}

A few comments are in order. By a C$^*$-correspondence over $A$ we mean a right Hilbert $A$-module together with a {\em nondegenerate} left action of $A$ on it. We have to explain why the left action on $(H_U\otimes B)^G$ is nondegenerate in the nonunital case. This is a consequence of the following simple lemma.

\begin{lemma}
If $\theta$ is a continuous left action of a compact quantum group $G$ on a C$^*$-algebra~$B$, then the fixed point algebra $A=B^G$ is a nondegenerate C$^*$-subalgebra of $B$.
\end{lemma}

\bp Let $\{e_s\}_s$ be an approximate unit in $A$. Define an $A$-valued inner product on $\B$ by $\langle x,y\rangle=E(x^*y)$. Then $xe_s\to x$ in the norm defined by this inner product for every $x\in\B$. By \cite[Lemma~2.5]{DCY0}, on every spectral subspace $B(U)\subset\B$ the norm defined by the inner product is equivalent to the C$^*$-norm. Therefore $xe_s\to x$ in the C$^*$-norm for every $x\in\B$, hence for every $x\in B$. 
\ep

C$^*$-correspondences over $A$ form a C$^*$-tensor category $\Hilb A$ with adjointable $A$-bilinear maps as morphisms. We emphasize that the isometries $F_{2,U,V}$ in the definition of the spectral functor are not claimed to be adjointable, and therefore formally they are not morphisms in $\Hilb A$. 

Finally, recall that two natural notions of isometry between Hilbert modules coincide: if~$M$ and~$N$ are right Hilbert $A$-modules, and $T\colon M\to N$ is an $A$-linear map such that $\|TX\|=\|X\|$ for all~$X\in M$, then $\langle TX,TY\rangle=\langle X,Y\rangle$ for all $X,Y\in M$, see e.g.~\cite[Theorem~3.5]{Lance}.

\bigskip

\section{Weak tensor functors}\label{s2}

Our goal is to give an abstract characterization of spectral functors. Here is the main definition.

\begin{definition}\label{dwutf}
Given a C$^*$-algebra $A$ and a strict C$^*$-tensor category $\CC$ with unit object $\un$, by a weak unitary tensor functor $\CC\to\Hilb A$ we mean a linear functor $F\colon\CC\to\Hilb A$ together with natural $A$-bilinear isometries $F_2=F_{2,U,V}\colon F(U)\otimes_A F(V)\to F(U\otimes V)$ such that the following conditions are satisfied:
\begin{itemize}
\item[(i)] $F(\un)=A$;
\item[(ii)] $F(T)^*=F(T^*)$ for any morphism $T$ in $\CC$;
\item[(iii)] $F_2\colon A\otimes_AF(U)\to F(\un\otimes U)=F(U)$ maps $a\otimes X$ into $aX$, and similarly $F_2\colon F(U)\otimes_A A\to F(U)$ maps $X\otimes a$ into $Xa$;
\item[(iv)] the diagrams
\begin{equation*}
\xymatrix{ F(U)\otimes_A F(V)\otimes_A F(W) \ar[d]_{\iota\otimes F_2}
\ar[r]^{\ \ F_2\otimes\iota}
 & F(U\otimes V)\otimes_A F(W)\ar[d]^{F_2}\\
F(U)\otimes_A F(V\otimes W) \ar[r]^{\ \ \ F_2}
& F(U\otimes V\otimes W)}
\end{equation*}
commute;
\item[(v)] for all objects $U$ and $V$ in $\CC$ and every vector $X\in F(U)$, the right $A$-linear map $S_X=S_{X,V}\colon F(V)\to F(U\otimes V)$ mapping $Y\in F(V)$ into $F_2(X\otimes Y)$ is adjointable, and the diagrams
\begin{equation*}
\xymatrix{ F(U\otimes V)\otimes_A F(W) \ar[d]_{S^*_X\otimes\iota}
\ar[r]^{\ \ \ F_2}
 & F(U\otimes V \otimes W)\ar[d]^{S^*_X}\\
F(V)\otimes_A F(W) \ar[r]^{\ \ \ F_2}
& F(V\otimes W)}
\end{equation*}
commute.
\end{itemize}
\end{definition}

Note that any unitary tensor functor $\CC\to\Hilb A$ defines a weak unitary tensor functor. In other words, if conditions (i)-(iv) are satisfied and the maps $F_2$ are surjective, then condition (v) is also satisfied. Indeed, the map $S_X\colon F(V)\to F(U\otimes V)$ is adjointable, because by assumption $F_2$ is unitary and the map $Y\mapsto X\otimes Y$ is adjointable, with adjoint given by $X'\otimes Y'\mapsto \langle X,X'\rangle Y'$. Since
$$
S_XF_2=F_2(S_X\otimes\iota)\colon F(V)\otimes_A F(W)\to F(U\otimes V\otimes W),
$$
by taking the adjoints we get $F_2^*S^*_X=(S^*_X\otimes\iota)F_2^*$. This is equivalent to commutativity of the diagram in (v) by unitarity of $F_2$.

Note also that if we consider $F$ simply as a functor into the category of vector spaces, then $S_X$ is a natural transformation from $F$ to $F(U\otimes\cdot)$, and so $S_X^*$ is a natural transformation from $F(U\otimes\cdot)$ to~$F$, or in other words,
\begin{equation} \label{enats}
S_X^*F(\iota\otimes T)=F(T)S^*_X
\end{equation}
for morphisms $T$ in $\CC$.

\medskip

Given a continuous left action of a compact quantum group $G$ on a C$^*$-algebra $B$ with fixed point algebra $A$, the associated spectral functor $\Rep G\to\Hilb A$ is a weak unitary tensor functor. Indeed, properties (i)-(iv) are immediate, while (v) follows by observing that the adjoint of the map
$$
S_X\colon (H_V\otimes B)^G\to (H_{U\times V}\otimes B)^G, \ \ Y\mapsto X_{13}Y_{23},
$$
is given by $S^*_XZ=X^*_{13}Z$. In other words, if $X=\sum_i\xi_i\otimes x_i$ and $Z=\sum_{i,j}\xi_i\otimes\zeta_j\otimes z_{ij}$ for orthonormal bases $\{\xi_i\}_i$ in $H_U$ and $\{\zeta_j\}_j$ in $H_V$, then
\begin{equation}\label{esxstar}
S^*_XZ=\sum_{i,j}\zeta_j\otimes x_i^*z_{ij}\in (H_V\otimes B)^G.
\end{equation}

The following is our main result.

\begin{theorem}\label{tmain}
Assume $G$ is a reduced compact quantum group and $A$ is a C$^*$-algebra. Then by associating to an action of $G$ on a C$^*$-algebra its spectral functor we get a bijection between isomorphism classes of triples $(B,\theta,\psi)$, where $\theta$ is a continuous left action of $G$ on a C$^*$-algebra $B$ and $\psi\colon A\to B$ is an embedding such that $B^G=\psi(A)$, and natural unitary monoidal isomorphism classes of weak unitary tensor functors $\Rep G\to\Hilb A$.
\end{theorem}

In the proof we will identify $A$ with $\psi(A)$ and simply talk about actions with fixed point algebra~$A$.

\medskip

The main part of the proof is, of course, a construction of an action from a weak unitary tensor functor $F\colon\Rep G\to\Hilb A$. We will define this action in a series of lemmas.

\medskip

Choose representatives $U_\alpha$ of isomorphism classes of irreducible unitary representations of $G$, and write~$H_\alpha$ instead of $H_{U_\alpha}$ for the underlying Hilbert spaces. We assume that there exists an index $e$ such that $U_e=\un$. Consider the space
$$
\B_F=\bigoplus_\alpha\bar H_\alpha\otimes F(U_\alpha).
$$
It will also be convenient to consider a much larger space. Choose a small C$^*$-tensor subcategory $\CC\subset\Rep G$ containing the objects $U_\alpha$, and then put
$$
\tilde\B_F=\bigoplus_U\bar H_U\otimes F(U),
$$
where the summation is over all objects in $\CC$. We have a canonical linear map $\pi\colon\tilde\B_F\to\B_F$ defined as follows. For a finite dimensional unitary representation $U$ of $G$, choose isometries $w_i\in\Mor(U_{\alpha_i},U)$ such that $\sum_iw_iw_i^*=\iota$. Then put
$$
\pi(\bar\xi\otimes X)=\sum_i\bar w^*_i\bar\xi\otimes F(w_i^*)X,
$$
where $\bar w_i\bar\zeta=\overline{w_i\zeta}$, so $\bar w^*_i\bar\xi=\overline{w_i^*\xi}$. This definition is independent of the choice of isometries $w_i$, since for any other choice $v_j$ there exists a unitary matrix $(u_{ij})_{i,j}$ such that $w_i=\sum_ju_{ij}v_j$. One property of~$\pi$ that we will regularly use, is that if $\bar\xi\otimes X\in\bar H_U\otimes F(U)$ and $w\in\Mor(U,V)$ is an isometry, then
\begin{equation} \label{epi}
\pi\big(\overline{w\xi}\otimes F(w)X\big)=\pi(\bar\xi\otimes X).
\end{equation}

Define a product on $\tilde\B_F$ by
$$
(\bar\xi\otimes X)\cdot(\bar\zeta\otimes Y)=\overline{(\xi\otimes\zeta)}\otimes F_2(X\otimes Y).
$$
It is immediate that this product is associative. Considering $\B_F$ as a subspace of $\tilde\B_F$, we define a product on $\B_F$ by
$$
xy=\pi(x\cdot y)\ \ \text{for}\ \ x,y\in\B_F.
$$

\begin{lemma}
The map $\pi\colon\tilde\B_F\to\B_F$ is a homomorphism, hence the product on $\B_F$ is associative.
\end{lemma}

\bp We have to check that $\pi(\pi(x)\cdot\pi(y))=\pi(x\cdot y)$ for all $x,y\in\tilde\B_F$. Take $x=\bar\xi\otimes X\in\bar H_U\otimes F(U)$, $y=\bar\zeta\otimes Y\in\bar H_V\otimes F(V)$ and choose isometries $u_i\in\Mor(U_{\alpha_i},U)$, $v_j\in\Mor(U_{\alpha_j},V)$ and $w_{ijk}\in\Mor(U_{\alpha_k},U_{\alpha_i}\times U_{\alpha_j})$ defining decompositions of $U$, $V$ and $U_{\alpha_i}\times U_{\alpha_j}$ into irreducibles. Then
$$
\pi(\pi(x)\cdot\pi(y))=\sum_{i,j,k}\overline{w_{ijk}^*(u_i^*\xi\otimes v_j^*\zeta)}\otimes F(w^*_{ijk})F_2(F(u_i^*)X\otimes F(v_j^*)Y),
$$
while
$$
\pi(x\cdot y)=\sum_{i,j,k}\overline{w_{ijk}^*(u_i^*\xi\otimes v_j^*\zeta)}\otimes F(w^*_{ijk}(u^*_i\otimes v_j^*))F_2(X\otimes Y).
$$
By naturality of $F_2$ these expressions are equal.
\ep

We can identify the space $\bar H_{e}\otimes F(U_e)=\bar\C\otimes A\subset\B_F$ with $A$. Under this identification, the space $A$, with its original product, becomes a subalgebra of $\B_F$. Furthermore, the left and right multiplications on $\bar H_\alpha\otimes F(U_\alpha)$ by elements of $A\subset \B_F$ are defined by the $A$-bimodule structure on~$F(U_\alpha)$, that is, the product on $\B_F$ has the property
$$
a(\bar\xi\otimes X)=\bar\xi\otimes aX,\ \ (\bar\xi\otimes X)a=\bar\xi\otimes Xa.
$$

Our next goal is to define an involution on $\B_F$. For a finite dimensional unitary representation $U$ of $G$ consider the standard solution $(R_U,\bar R_U)$ of the conjugate equations for $U$ defined by~\eqref{econj}.

\begin{lemma}\label{lbullet}
For every $X\in F(U)$ there exists a unique element $X^\bullet\in F(\bar U)$ such that
$$
\langle X^\bullet,Y\rangle=F(\bar R_U^*)F_2(X\otimes Y)\ \ \text{for all}\ \ Y\in F(\bar U).
$$
If the C$^*$-algebra $A$ is unital, then $X^\bullet=S^*_XF(\bar R_U)(1)$. We also have
$$
\langle X,Y\rangle=F(R_U^*)F_2(X^\bullet\otimes Y)\ \ \text{for all}\ \ Y\in F(U).
$$
\end{lemma}

\bp The uniqueness is clear. In order to prove the existence assume first that $A$ is unital. We then have
$$
\langle S^*_XF(\bar R_U)(1),Y\rangle=\langle F(\bar R_U)(1),F_2(X\otimes Y)\rangle
=F(\bar R_U^*)F_2(X\otimes Y),
$$
so $X^\bullet=S^*_XF(\bar R_U)(1)$. If $A$ is nonunital, then a similar computation shows that for any $X\in F(U)$ and $a\in A$ the element $(aX)^\bullet$ exists and $(aX)^\bullet=S^*_XF(\bar R_U)(a^*)$. But this is enough, since by Cohen's factorization theorem any element of $F(U)$ has the form $aX$.

To prove the last statement in the formulation, assume once again that $A$ is unital, the nonunital case requires only a minor modification. For $Y\in F(U)$ we compute:
\begin{align*}
F(R^*_U)F_2(X^\bullet\otimes Y)
&=F(R^*_U)F_2(S^*_XF(\bar R_U)(1)\otimes Y)\\
&=F(R^*_U)S^*_XF_2(F(\bar R_U)(1)\otimes Y)\\
&=S^*_XF(\iota\otimes R^*_U)F(\bar R_U\otimes\iota)Y&\text{by \eqref{enats}}\\
&=S^*_XY.
\end{align*}
Here $S_X$ is the map $A\to F(U)$, $a\mapsto F_2(a\otimes X)=aX$, so $S^*_XY=\langle X,Y\rangle$.
\ep

This lemma implies that the correspondences $F(U)$ and $F(\bar U)$ are, in some sense, dual to each other. In general, this is not the duality in the C$^*$-categorical sense. Already the simplest examples, such as the spectral functor associated with the action of $\T$ by rotations on the unit disk, show that the objects $F(U)$ do not necessarily have conjugates in $\Hilb A$.

Similarly, for every vector $\xi\in H_U$ define a vector $\xi^\bullet\in H_{\bar U}$ by
$$
\xi^\bullet=(\iota\otimes\bar\xi)R_U(1)=\overline{\rho^{-1/2}\xi},\ \text{so}\ \ (\zeta,\xi^\bullet)=R^*_U(\zeta\otimes\xi)\ \ \text{for all}\ \ \zeta\in H_{\bar U}.
$$
Define an anti-linear map $\tilde\B_F\to\tilde\B_F$, $x\mapsto x^\bullet$, by
$$
(\bar\xi\otimes X)^\bullet=\overline{\xi^\bullet}\otimes X^\bullet.
$$
For $x\in\B_F$ put $x^*=\pi(x^\bullet)$. On $A\subset\B_F$ this clearly coincides with the involution on $A$. Although we will not need this, we remark that it is not difficult to show that the particular choice of solutions $(R_U,\bar R_U)$ was not important for defining the involution on $\B_F$, in the sense that for every $x\in\tilde \B_F$ the element $\pi(x^\bullet)$ is independent of any choices.

\begin{lemma}
The map $x\mapsto x^*$ defines an involution on the algebra $\B_F$, and for every $x\in\tilde\B_F$ we have $\pi(x)^*=\pi(x^\bullet)$.
\end{lemma}

\bp We start by proving the second part. We have to show that $\pi(\pi(x)^\bullet)=\pi(x^\bullet)$. Take an element $x=\bar\xi\otimes X\in\bar H_U\otimes F(U)$. Choose isometries $w_i\in\Mor(U_{\alpha_i},U)$ defining a decomposition of $U$ into irreducibles. Write $R_i$ for $R_{U_{\alpha_i}}$ and $\bar R_i$ for $\bar R_{U_{\alpha_i}}$. Then $R_U=\sum_i(\bar w_i\otimes w_i)R_i$ and $\bar R_U=\sum_i(w_i\otimes \bar w_i)\bar R_i$. For any $Y\in F(\bar U)$ we have
$$
F(\bar R_U^*)F_2(X\otimes Y)=\sum_iF(\bar R_i^*)F_2(F(w_i^*)X\otimes F(\bar w_i^*) Y)
=\sum_i\langle (F(w_i^*)X)^\bullet,F(\bar w_i^*) Y\rangle,
$$
so $X^\bullet=\sum_i F(\bar w_i)(F(w_i^*)X)^\bullet$. We also have $\xi^\bullet=\sum_i \bar w_i(w_i^*\xi)^\bullet$. Therefore
$$
x^\bullet=\sum_i \bar{\bar w}_i\overline{(w_i^*\xi)^\bullet}\otimes F(\bar w_i)(F(w_i^*)X)^\bullet.
$$
Applying $\pi$ and using \eqref{epi} we get
$$
\pi(x^\bullet)=\sum_i\pi\big(\overline{(w_i^*\xi)^\bullet}\otimes(F(w_i^*)X)^\bullet\big)
=\pi(\pi(x)^\bullet).
$$


We next prove anti-multiplicativity of the map $*$ on $\B_F$. For this it suffices to check that for all $x,y\in\tilde\B_F$ we have $\pi((x\cdot y)^\bullet)=\pi(y^\bullet\cdot x^\bullet)$. Take $x=\bar\xi\otimes X\in\bar H_U\otimes F(U)$ and $y=\bar\zeta\otimes Y\in\bar H_V\otimes F(V)$. The unitary $\sigma\colon H_{\bar V}\otimes H_{\bar U}\to H_{\overline{U\times V}}$ mapping $\bar\eta\otimes\bar\vartheta$ into $\overline{\vartheta\otimes\eta}$ defines an equivalence between $\bar V\times\bar U$ and $\overline{U\times V}$, and we have
$$
R_{U\times V}=(\sigma\otimes\iota\otimes\iota)(\iota\otimes R_U\otimes\iota)R_V\ \ \text{and}\ \ \bar R_{U\times V}=(\iota\otimes\iota\otimes\sigma)(\iota\otimes \bar R_V\otimes\iota)\bar R_U.
$$
Assuming that $A$ is unital we compute:
\begin{align*}
F_2(X\otimes Y)^\bullet &=S_{F_2(X\otimes Y)}^*F(\bar R_{U\times V})(1)&\text{by Lemma~\ref{lbullet}}\\
&=S^*_YS_X^*F(\iota\otimes\iota\otimes\sigma)F(\iota\otimes \bar R_V\otimes\iota)F(\bar R_U)(1)&\text{as\ \ }S_{F_2(X\otimes Y)}=S_XS_Y \\
&=F(\sigma)S^*_YF(\bar R_V\otimes\iota)S^*_XF(\bar R_U)(1)&\text{by \eqref{enats}}\\
&=F(\sigma)S^*_YF(\bar R_V\otimes\iota)(X^\bullet)\\
&=F(\sigma)S^*_YF_2(F(\bar R_V)(1)\otimes X^\bullet)\\
&=F(\sigma)F_2(S^*_YF(\bar R_V)(1)\otimes X^\bullet )\\
&=F(\sigma)F_2(Y^\bullet\otimes X^\bullet).
\end{align*}
In the nonunital case we get the same identity by replacing $X$ and $Y$ by elements of the form $aX$ and $bY$, see the proof of Lemma~\ref{lbullet}. We also have $(\xi\otimes\zeta)^\bullet=\sigma(\zeta^\bullet\otimes\xi^\bullet)$. Therefore
$$
(x\cdot y)^\bullet=(\bar\sigma\otimes F(\sigma))\big(\overline{(\zeta^\bullet\otimes\xi^\bullet)}\otimes F_2(Y^\bullet\otimes X^\bullet)\big)=(\bar\sigma\otimes F(\sigma))(y^\bullet\cdot x^\bullet).
$$
Applying $\pi$ we get $\pi((x\cdot y)^\bullet)=\pi(y^\bullet\cdot x^\bullet)$.

\smallskip

It remains to show that the map $x\mapsto x^*$ on $\B_F$ is involutive. Equivalently, we have to show that $\pi(x^{\bullet\bullet})=\pi(x)$ for all $x\in\tilde\B_F$. Take an element $x=\bar\xi\otimes X\in\bar H_U\otimes F(U)$. Consider the unitary $u\colon H_U\to H_{\bar{\bar U}}$ mapping $\zeta$ into $\bar{\bar\zeta}$. Then $\bar R_{\bar U}=(\iota\otimes u)R_U$. Hence, for any $Y\in F(\bar{\bar U})$,
$$
\langle X^{\bullet\bullet},Y\rangle=F(\bar R_{\bar U}^*)F_2(X^\bullet\otimes Y)
=F(R^*_U)F_2(X^\bullet\otimes F(u^*)Y)=\langle X,F(u^*)Y\rangle,
$$
where the last equality follows from Lemma~\ref{lbullet}. Thus $X^{\bullet\bullet}=F(u)X$. We also have $\xi^{\bullet\bullet}=\bar{\bar\xi}=u\xi$. Therefore
$$
x^{\bullet\bullet}=(\bar u\otimes F(u))x,
$$
and applying $\pi$ we get $\pi(x^{\bullet\bullet})=\pi(x)$.
\ep

We next define a linear map $\theta_F\colon\B_F\to\C[G]\otimes\B_F$ by
$$
\theta_F(\bar\xi\otimes X)=(U^c_\alpha)^*_{21}(1\otimes\bar\xi\otimes X)\ \ \text{for}\ \ \bar\xi\otimes X\in\bar H_\alpha\otimes F(U_\alpha).
$$
In other words, if we fix an orthonormal basis $\{\xi_i\}_i$ in $H_\alpha$ and write $U_\alpha$ as a matrix $(u_{ij})_{i,j}$, then
$$
\theta_F(\bar\xi_i\otimes X)=\sum_ju_{ij}\otimes\bar\xi_j\otimes X.
$$

\begin{lemma}
The map $\theta_F$ defines a left algebraic action of $G$ on $\B_F$ with fixed point algebra $A$.
\end{lemma}

\bp Clearly, the map $\theta_F$ turns $\B_F$ into a comodule over $(\C[G],\Delta)$ with fixed point subcomodule~$A$. 

\smallskip

In order to show that $\theta_F$ is a homomorphism, observe first that we have a left comodule structure $\tilde\theta_F\colon\tilde\B_F\to \C[G]\otimes\tilde\B_F$ on $\tilde\B_F$ defined in the same way as for $\B_F$, so $\theta_F(\bar\xi\otimes X)=(U^c)^*_{21}(1\otimes\bar\xi\otimes X)$ for $\bar\xi\otimes X\in\bar H_U\otimes F(U)$. Then $\pi\colon\tilde\B\to\B_F$ is a comodule map, since if $w\in\Mor(U,V)$, then $U^{c*}(\bar w\otimes1)=(\bar w\otimes1)V^{c*}$. Using that $(U\times V)^{c*}=(U^{c})^*_{13}(V^{c})^*_{23}$, modulo identification of $\overline{H_U\otimes H_V}$ with $\bar H_U\otimes\bar H_V$, it is easy to see that $\tilde\theta_F$ is a homomorphism. Hence $\theta_F$ is also a homomorphism.

\smallskip

Next let us check that $\theta_F$ is $*$-preserving. It suffices to show that $\tilde\theta_F(x)^{*\otimes\bullet}=\tilde\theta_F(x^\bullet)$ for $x\in \bar H_U\otimes F(U)\subset\tilde\B_F$. Fixing an orthonormal basis $\{\xi_i\}_i$ in $H_U$ and identifying $\bar{\bar H}_U$ with $H_U$, we get
$$
\tilde\theta_F(\bar\xi_i\otimes X)^{*\otimes\bullet}=\sum_ju^*_{ij}\otimes\rho^{-1/2}\xi_j\otimes X^\bullet\ \ \text{and}\ \ \tilde\theta_F((\bar\xi_i\otimes X)^\bullet)=(\bar U^c)^{*}_{21}(1\otimes\rho^{-1/2}\xi_i\otimes X^\bullet).
$$
Since $\bar U^{c*}=(\rho^{-1/2}\otimes1)U^*(\rho^{1/2}\otimes1)$, these expressions coincide.

\smallskip

It remains to show that $\B_F$ is a right pre-Hilbert $A$-module with inner product $\langle x,y\rangle=E(x^*y)$, where  $E=(h\otimes\iota)\theta_F$, and the left action of $A$ on $\B_F$ by multiplication is bounded. This will follow immediately, if we can show that the spaces $\bar H_\alpha\otimes F(U_\alpha)$ are mutually orthogonal and
$$
\langle\bar\xi\otimes X,\bar\zeta\otimes Y\rangle=\frac{1}{\dim_qU_\alpha}(\rho^{-1}\xi,\zeta)\langle X,Y \rangle\ \ \text{for}\ \ \xi,\zeta\in H_\alpha,\ X,Y\in F(U_\alpha).
$$
Note that if $z=\bar\eta\otimes Z\in\bar H_U\otimes F(U)$, then $E(\pi(z))=\sum_i\overline{w^*_i\eta}\otimes F(w^*_i)Z$, where $w_i\in\Mor(\un,U)$ are isometries such that $\sum_iw_iw_i^*$ is the projection onto the isotypic component of $U$ corresponding to the trivial representation. This clearly implies mutual orthogonality of the spaces $\bar H_\alpha\otimes F(U_\alpha)$. If $U=\bar U_\alpha\times U_\alpha$, then the only isometry in $\Mor(\un,U)$, up to a phase factor, is $(\dim_qU_\alpha)^{-1/2}R_\alpha$, where $R_\alpha=R_{U_\alpha}$. Therefore for $\xi,\zeta\in H_\alpha$ and $X,Y\in F(U_\alpha)$ we have
$$
\langle\bar\xi\otimes X,\bar\zeta\otimes Y\rangle=\frac{1}{\dim_qU_\alpha}\overline{R^*_\alpha(\xi^\bullet\otimes\zeta)}F(R^*_\alpha)F_2(X^\bullet\otimes Y).
$$
By Lemma~\ref{lbullet} we have $F(R^*_\alpha)F_2(X^\bullet\otimes Y)=\langle X,Y\rangle$. Using that $\xi^\bullet=\overline{\rho^{-1/2}\xi}$ it is also straightforward to check that $R^*_\alpha(\xi^\bullet\otimes\zeta)=(\zeta,\rho^{-1}\xi)$. This finishes the proof of the lemma.
\ep

As we discussed in the previous section, an algebraic action of $G$ on $\B$ uniquely defines a completion~$B_F$ of~$\B_F$ carrying a continuous action of the reduced form of $G$. Therefore the previous lemma finishes our construction of a continuous action from a weak unitary tensor functor.

\bp[Proof of Theorem~\ref{tmain}]
It is clear that isomorphic actions produce  naturally unitarily monoidally  isomorphic weak unitary tensor functors, and naturally unitarily  monoidally isomorphic weak unitary tensor functors produce isomorphic actions. It remains to show that up to isomorphisms the constructions are inverse to each other.

\smallskip

Assume $\theta$ is a continuous left action of $G$ on a C$^*$-algebra $B$ with fixed point algebra $A$. Let $F$ be the associated spectral functor and $\B\subset B$ be the subalgebra of regular elements. Consider the algebraic action $\theta_F$ of $G$ on $\B_F$ defined by $F$ as described above. We have a linear isomorphism
$$
\B_F\cong\B\ \ \text{mapping}\ \ \pi(\bar\xi\otimes X)\in\B_F\ \ \text{into}\ \ (\bar\xi\otimes\iota)(X)\in\B
$$
for $\bar\xi\otimes X\in\bar H_U\otimes (H_U\otimes B)^G$.
It is easy to see that this is a $G$-equivariant isomorphism of algebras. It is a bit less obvious that this isomorphism is $*$-preserving. In order to show this, fix an irreducible representation $U_\alpha$ and an orthonormal basis $\{\xi_i\}_i$ in $H_\alpha$. Consider an element $X=\sum_i\xi_i\otimes x_i\in (H_\alpha\otimes B)^G$. Writing $\bar R_\alpha$ for $\bar R_{U_\alpha}$, assuming for simplicity that $A$ is unital and using Lemma~\ref{lbullet} and identity~\eqref{esxstar} for $S_X^*$, we get
$$
X^\bullet=S^*_XF(\bar R_\alpha)(1)=S^*_X\left(\sum_j\rho^{1/2}\xi_j\otimes\bar\xi_j\otimes1\right)
=\sum_{i,j}(\rho^{1/2}\xi_j,\xi_i)\bar\xi_j\otimes x^*_i=\sum_i\overline{\rho^{1/2}\xi_i}\otimes x_i^*.
$$
From this we see that the image of the element $(\bar\xi\otimes X)^*=\pi\Big(\overline{\overline{\rho^{-1/2}\xi}}\otimes X^\bullet\Big)\in\B_F$ in $\B$ equals
$$
\sum_i(\rho^{-1/2}\xi,\rho^{1/2}\xi_i)x_i^*=\left(\sum_i(\xi_i,\xi)x_i\right)^*,
$$
so the isomorphism $\B_F\cong\B$ is indeed $*$-preserving.

\smallskip

Now conversely, assume we start with a weak unitary tensor functor $F$, consider the action $\theta_F$ of $G$ on $B_F$, and define the corresponding spectral functor $F'$. It is easy to see that if we fix an irreducible representation $U_\alpha$ and an orthonormal basis $\{\xi_i\}_i$ in $H_\alpha$, then the dense subspace $(H_\alpha\otimes\B_F)^G$ of $F'(U_\alpha)=(H_\alpha\otimes B_F)^G$ consists of vectors of the form $\sum_i(\xi_i\otimes\bar\xi_i\otimes X)$, with $X\in F(U_\alpha)$. We have the obvious $A$-bilinear map $F(U_\alpha)\to F'(U_\alpha)$ with dense image, mapping $X$ into $\sum_i(\xi_i\otimes\bar\xi_i\otimes X)$. Let us check that this map is isometric. Taking vectors $X'=\sum_i(\xi_i\otimes\bar\xi_i\otimes X)$ and $Y'=\sum_i(\xi_i\otimes\bar\xi_i\otimes Y)$ in $F'(U_\alpha)$, and writing $R_\alpha$ for $R_{U_\alpha}$, we compute:
$$
\langle X',Y'\rangle=\sum_i(\bar\xi_i\otimes X)^*(\bar\xi_i\otimes Y)
=\pi\left(\sum_i(\overline{\overline{\rho^{-1/2}\xi_i}}\otimes X^\bullet)\cdot(\bar\xi_i\otimes Y)\right)
=\pi\big(\overline{R_\alpha(1)}\otimes F_2(X^\bullet\otimes Y)\big).
$$
Since, up to a scalar factor, $R_\alpha$ is an isometry in $\Mor(\un,\bar U_\alpha\times U_\alpha)$, the last expression equals
$$
F(R^*_\alpha)F_2(X^\bullet\otimes Y)=\langle X,Y\rangle
$$
by Lemma~\ref{lbullet}. Thus we get unitary isomorphisms $F(U_\alpha)\cong F'(U_\alpha)$. These isomorphisms for all $\alpha$ extend uniquely to a natural unitary isomorphism between the functors $F$ and $F'$. It is straightforward to check that this isomorphism is monoidal.
\ep

\bigskip

\section{Module categories}

In this section we will give a different categorical description of actions in terms of module categories. Recall that given a C$^*$-tensor category $\CC$, a right $\CC$-module C$^*$-category is a C$^*$-category~$\M$ equipped with a bilinear unitary functor $\otimes\colon \M\times\CC\to\M$ together with natural unitary isomorphisms $\phi\colon (M\otimes U)\otimes V\to M\otimes (U\otimes V)$ and $e\colon M\otimes\un\to M$ satisfying certain coherence relations, see~\cite{DCY} for details. If $\CC$ is strict, a module category $\M$ is called strict if $\phi$ and $e$ are the identity morphisms. Any module category over a strict C$^*$-tensor category is equivalent to a strict one. 

In the following discussion we will tacitly assume that the C$^*$-categories that we consider have subobjects, meaning that for every projection $p$ in $\End(M)$ there exists an object $N$ and an isometry $v\in \Mor(N,M)$ such that $vv^*=p$. This is a very mild assumption, as we can always complete a C$^*$-category with respect to subobjects.

\smallskip

Assume we are given a continuous left action of a reduced compact quantum group $G$ on a {\em unital} C$^*$-algeb\-ra~$B$. Following~\cite{DCY}, consider the category $\D_B$ of unitary $G$-equivariant finitely generated right Hilbert $B$-modules. By definition, the morphisms in $\D_B$ are $G$-equivariant maps of Hilbert $B$-modules. Since we consider only finitely generated Hilbert modules, such maps are automatically adjointable, so $\D_B$ is a C$^*$-category. It is a strict right $(\Rep G)$-module C$^*$-category: given a right Hilbert $B$-module $M$ with the action of~$G$ given by an isometry $\delta_M\colon M\to C(G)\otimes M$, we define $M\otimes U$ as the Hilbert $B$-module $M\otimes H_U$ with the action of $G$ given by $x\otimes \xi\mapsto U^*_{31}(\delta_M(x)\otimes\xi)$. Note that for $M=B$ the module $B\otimes U$ is, up to identification of $H_U\otimes B$ with $B\otimes H_U$, the same equivariant module $H_U\otimes B$ that we considered in Section~\ref{s1}. The module $B$ generates the category~$\D_B$, in the sense that any object in $\D_B$ is a subobject of $B\otimes U$ for some $U$. In other words, any $G$-equivariant finitely generated right Hilbert $B$-module $M$ is isomorphic to a direct summand of $B\otimes H_U$ for some $U$, see \cite[Section~3.2]{Ver} or \cite[Lemma~3.2]{NTK}.

Let $F$ be the spectral functor associated with the action of $G$ on $B$. We have canonical isomorphisms
$$
F(U)=(H_U\otimes B)^G\cong\Mor(B,B\otimes U)
$$
that map $\sum_i(\xi_i\otimes x_i)\in(H_U\otimes B)^G$ into the morphism $x\mapsto \sum_i(x_ix\otimes\xi_i)$. For ergodic actions, these isomorphisms, modulo some identifications in terms of Frobenius reciprocity, were already used in \cite[Section~6]{DCY} to identify algebras constructed in~\cite{DCY} with those defined by Pinzari and Roberts~\cite{PR}. In other words, the key relation between the results in \cite{DCY} and \cite{PR} can be described by saying that given a $(\Rep G)$-module C$^*$-category $\M$ and a simple object $M$ in $\M$, the functor $U\mapsto\Mor(M,M\otimes U)$ has all the properties of a spectral functor.  With our characterization of spectral functors this becomes almost immediate. Specifically, and more generally, we have the following.

\begin{proposition}
Assume $\M$ is a strict right module C$^*$-category over a strict C$^*$-tensor category~$\CC$. Take an object $M\in\M$ and consider the unital C$^*$-algebra $A=\End(M)$. Then the following defines a weak unitary tensor functor $\CC\to\Hilb A$:
$$
F(U)=\Mor(M,M\otimes U),\ \
$$
with the right $A$-module structure on $F(U)$ given by composition of morphisms, the left $A$-module structure by $aX=(a\otimes\iota)X$ and the inner product by $\langle X,Y\rangle=X^*Y$, the action of $F$ on morphisms is defined by $F(T)X=(\iota\otimes T)X$, and
$$
F_2\colon F(U)\otimes_A F(V)\to F(U\otimes V)
$$
is given by $X\otimes Y\mapsto (X\otimes\iota)Y$.
\end{proposition}

\bp This is a routine verification. We only remark that the adjoint of the map
$$
S_X\colon F(V)\to F(U\otimes V),\ \ Y\mapsto (X\otimes\iota)Y,
$$
is obviously given by $S^*_XZ=(X^*\otimes\iota)Z$.
\ep

If the object $M$ happens to be generating, then we can reconstruct the whole category $\M$ from the functor $F$. For $\CC=\Rep G$ this gives the following result.

\begin{proposition}
Assume $G$ is a reduced compact quantum group and $\M$ is a strict right $(\Rep G)$-module C$^*$-category generated by an object $M$. Put $A=\End(M)$ and consider the weak unitary tensor functor $F\colon\Rep G\to\Hilb A$ defined by the object $M$ as described in the previous proposition. Let $\theta\colon B\to C(G)\otimes B$ be the continuous action corresponding to this functor by Theorem~\ref{tmain}. Then~$\D_B$ is unitarily equivalent, as a $(\Rep G)$-module C$^*$-category, to $\M$, via an equivalence that maps the generator $B\in\D_B$ into $M$.
\end{proposition}

\bp Consider the functor $F'\colon\D_B\to\Hilb A$ defined by the object $B\in\D_B$. By the above discussion, it is naturally unitarily monoidally   isomorphic to the spectral functor associated with the action of $G$ on $B$, hence to $F$. Let $\psi\colon F'\to F$ be such an isomorphism. Note that we automatically have that $\psi\colon A=F'(\un)\to F(\un)=A$ is the identity map, since it is a bimodule map such that $\psi F_2=F_2'(\psi\otimes\psi)$.

\smallskip

Consider the full subcategories $\tilde\D_B\subset\D_B$ and $\tilde\M\subset\M$ consisting of objects $B\otimes U$ and $M\otimes U$, respectively. We want to define a functor $E\colon\tilde\D_B\to\tilde\M$. On objects we put $E(B\otimes U)=M\otimes U$. For morphisms $T\in\Mor(B,B\otimes U)$ we put $E(T)=\psi(T)$. More generally, given two finite dimensional unitary representations $U$ and $V$, we have Frobenius reciprocity isomorphisms
$$
\Mor(B\otimes U,B\otimes V)\to\Mor(B,B\otimes V\otimes\bar U),\ \ T\mapsto(T\otimes\iota)(\iota\otimes\bar R_U),
$$
with inverse $S\mapsto(\iota\otimes\iota\otimes R^*_U)(S\otimes\iota)$. We also have similar isomorphisms in $\tilde\M$. Hence we can define linear isomorphisms
$$
E\colon\Mor(B\otimes U,B\otimes V)\to\Mor(M\otimes U,M\otimes V)
$$
by $E(T)=(\iota\otimes\iota\otimes R^*_U)\big(\psi\big((T\otimes\iota)(\iota\otimes\bar R_U)\big)\otimes\iota\big)$.

Before we turn to the proof that $E$ is indeed a functor, let us make two observations. The first one is that given a morphism $T\colon B\otimes U\to B\otimes V$ and a finite dimensional unitary representation~$W$ of~$G$, for the morphism $T\otimes\iota_W=T\otimes\iota\colon B\otimes U\otimes W\to B\otimes V\otimes W$ we have
\begin{equation} \label{ee1}
E(T\otimes\iota)=E(T)\otimes\iota.
\end{equation}
The second observation is that given a morphism $T\colon B\otimes U\to B\otimes V$ and a morphism $S\colon V\to W$, we have
\begin{equation} \label{ee2}
E((\iota\otimes S)T)=(\iota\otimes S)E(T).
\end{equation}
Both claims follow easily from naturality of $\psi$, which means that for any morphisms $T\colon B\to B\otimes U$ and $S\colon U\to V$ we have $\psi((\iota\otimes S)T)=(\iota\otimes S)\psi(T)$.

Consider now morphisms $R\colon B\otimes U\to B\otimes V$ and $T\colon B\otimes V\to B\otimes W$, and define the morphisms $P=(R\otimes\iota)(\iota\otimes\bar R_U)\colon B\to B\otimes V\otimes\bar U$ and $S=(T\otimes\iota)(\iota\otimes\bar R_V)\colon B\to B\otimes W\otimes\bar V$. We then have
\begin{align*}
TR&=(\iota_B\otimes\iota_W\otimes R^*_V)(S\otimes\iota_V)(\iota_B\otimes\iota_V\otimes R^*_U)(P\otimes\iota_U)\\
&=(\iota_B\otimes\iota_W\otimes R^*_V)(\iota_B\otimes\iota_W\otimes\iota_{\bar V}\otimes\iota_V\otimes R^*_U)( (S\otimes\iota\otimes\iota)P\otimes\iota_U)\\
&=(\iota_B\otimes\iota_W\otimes R^*_V\otimes R^*_U)(F'_2(S\otimes P)\otimes\iota_U),
\end{align*}
where $F'_2(S\otimes P)=(S\otimes\iota\otimes\iota)P\colon B\to B\otimes W\otimes\bar V\otimes V\otimes\bar U$. A similar computation gives
$$
E(T)E(R)=(\iota_M\otimes\iota_W\otimes R^*_V\otimes R^*_U)(F_2(\psi(S)\otimes \psi(P))\otimes\iota_U).
$$
From this we immediately get that $E(TR)=E(T)E(R)$ using \eqref{ee2}, \eqref{ee1} and monoidality of $\psi$, which means that $\psi F'_2(S\otimes P)=F_2(\psi(S)\otimes \psi(P))$. Therefore $E$ is a functor. Since it is surjective on objects and fully faithful, it is an equivalence of the linear categories $\tilde\D_\B$ and $\tilde\M$.

\smallskip

Let us show next that the equivalence $E$ is unitary, that is, $E(T)^*=E(T^*)$ on morphisms. Let us check this first for $T\in\Mor(B,B\otimes U)$. Since $\psi$ is unitary, for any $S\in\Mor(B,B\otimes U)$ we have
$$
E(T)^*E(S)=\psi(T)^*\psi(S)=\langle\psi(T),\psi(S)\rangle=\langle T,S\rangle=T^*S=E(T^*S)=E(T^*)E(S).
$$
Since this is true for all $S$, we conclude that $E(T)^*=E(T)^*$. By virtue of \eqref{ee1} we then also get $E(T\otimes\iota)^*=E((T\otimes\iota)^*)$. But any morphism in $\tilde\D_B$ is a composition of such a morphism $T\otimes\iota_W$ and a morphism of the form $\iota_M\otimes S$ for some morphism $S$ in $\Rep G$. Since as a particular case of~\eqref{ee2} we have
$
E(\iota\otimes S)^*=\iota\otimes S^*=E((\iota\otimes S)^*),
$
it follows that $E$ is unitary.

\smallskip

Next, from \eqref{ee1} and \eqref{ee2} we see that if we define
$$
E_2=E_{2,B\otimes U,V}\colon E(B\otimes U)\otimes V\to E(B\otimes U\otimes V)
$$
to be the identity maps, then we get a natural isomorphism of bilinear functors $E(\cdot)\otimes\cdot$ and~$E(\cdot\otimes\cdot)$. Therefore the pair $(E,E_2)$ defines a unitary equivalence of $(\Rep G)$-module categories $\tilde\D_B$ and $\tilde\M$.

Finally, since $\D_B$ and $\M$ are completions of these categories with respect to subobjects, the equivalence between $\tilde\D_B$ and $\tilde\M$ extends uniquely, up to a natural unitary isomorphism, to a unitary equivalence between the $(\Rep G)$-module C$^*$-categories $\D_B$ and $\M$.
\ep

This leads to the main theorem of this section, a generalization of results of De Commer and Yamashita~\cite{DCY} to the nonsemisimple/nonergodic case.

\begin{theorem}\label{tmain2}
Assume $G$ is a reduced compact quantum group. Then by associating to an action of~$G$ on a unital C$^*$-algebra $B$ the $(\Rep G)$-module category $\D_B$ with generator $B$, we get a bijection between isomorphism classes of continuous left actions of $G$ on unital C$^*$-algebras and unitary equivalence classes of pairs $(\M,M)$, where $\M$ is a right $(\Rep G)$-module C$^*$-category and $M$ is a generating object in $\M$.
\end{theorem}

\bp In view of the above proposition we only have to show that two actions of $G$ on unital C$^*$-algebras $B$ and $C$ are isomorphic if and only if the pairs $(\D_B,B)$ and $(\D_C,C)$ are unitarily equivalent. Given such an equivalence, we first of all get an isomorphism $B^G=\End_{\D_B}(B)\cong\End_{\D_C}(C)=C^G$. Modulo the identification of $B^G$ with $C^G$ using this isomorphism, we then also get a natural unitary isomorphism between the spectral functors associated with our actions. Hence the actions are isomorphic by the easy part of Theorem~\ref{tmain}. Conversely, it is clear that isomorphic actions produce unitarily equivalent pairs.
\ep

As in \cite{DCY}, this result can also be formulated in terms of Morita equivalent actions.

\begin{corollary}
For any reduced compact quantum group $G$, there is a bijection between Morita equivalence classes of continuous left actions of $G$ on unital C$^*$-algebras and unitary equivalence classes of singly generated right $(\Rep G)$-module C$^*$-categories.
\end{corollary}

\bp It suffices to show that two actions of $G$ on unital C$^*$-algebras $B$ and $C$ are Morita equivalent if and only if the $(\Rep G)$-module C$^*$-categories $\D_B$ and $\D_C$ are unitarily equivalent. In one direction this is obvious:
if a $C$-$B$-bimodule $M$ defines the Morita equivalence, then $\D_C$ and $\D_B$ are unitarily equivalent, via an equivalence that maps $N\in\D_C$ into $N\otimes_C M$. Conversely, assume we have a unitary equivalence $E\colon\D_C\to\D_B$ of $(\Rep G)$-module categories. Consider the right Hilbert $B$-module $M=E(C)$. Since it is a generating object in $\D_B$, and therefore the right Hilbert $B$-module~$B$ can be isometrically embedded into $M\otimes H_U$ for some representation $U$, the module $M$ must be full. Thus the action of $G$ on $B$ is Morita equivalent to the action of $G$ on $C'=\End_B(M)$. The $(\Rep G)$-module C$^*$-categories $\D_B$ and $\D_{C'}$ are unitarily equivalent, via an equivalence that maps $M\in\D_B$ into $C'\in\D_{C'}$. It follows that the $(\Rep G)$-module C$^*$-categories~$\D_C$ and~$\D_{C'}$ are unitarily equivalent, via an equivalence that maps~$C$ into~$C'$. By Theorem~\ref{tmain2} this implies that the actions of $G$ on $C$ and $C'$ are isomorphic, so the actions of $G$ on $B$ and $C$ are Morita equivalent.
\ep

\begin{remark}
In the above proof we used that if an equivariant right Hilbert $B$-module $M$ is a generating object in $\D_B$, then it is full. The proof implies that the converse is also true, since if~$M$ is full, then $M$ is the image of the generating object $C$ in $\D_C$, where $C=\End_B(M)$, under the equivalence of categories $\D_C$ and $\D_B$ defined by $M$, hence $M$ is a generating object in $\D_B$. Somewhat more explicitly this can also be proved as follows. It suffices to show that $B\in\D_B$ is a subobject of $M\otimes V$ for some finite dimensional unitary representation~$V$ of~$G$. Replacing $M$ by $M^n$ we may assume that there exists a vector $X\in M$ such that the element $\langle X,X\rangle\in B$ is invertible. Furthermore, since the union of spectral subspaces of~$M$ is dense in $M$, we may assume that $X$ lies in a spectral subspace of $M$ corresponding to some representation~$U$. In other words, there exist an orthonormal basis  $\{\xi_i\}_i$ in $H_U$ and vectors $X_i\in M$ such that $\delta_M(X_i)=\sum_ju^*_{ij}\otimes X_j$ and one of the inner products  $\langle X_i,X_i\rangle\in B$ is invertible. Consider the vector
$$
Y=\sum_iX_i\otimes\overline{\rho^{1/2}\xi_i}\in M\otimes H_{\bar U}.
$$
Then $Y$ is invariant and
$\langle Y,Y\rangle=\sum_{i,j}\langle X_i,X_j\rangle(\rho\xi_i,\xi_j)$. Since the matrix $((\rho\xi_i,\xi_j))_{i,j}$ is positive and invertible, and the matrix $(\langle X_i,X_j\rangle)_{i,j}$ is positive, there exists a constant $c>0$ such that
$$
\langle Y,Y\rangle\ge c\sum_i\langle X_i,X_i\rangle.
$$
It follows that $\langle Y,Y\rangle$ is invertible, so the map $B\ni x\mapsto Y\langle Y,Y\rangle^{-1/2}x$ gives an equivariant isometric embedding of $B$ into $M\otimes H_{\bar U}$.
\end{remark}

\bigskip

\section{Examples and applications}

The data provided by a weak unitary tensor functor and the construction of the corresponding C$^*$-algebra are reminiscent of various crossed product type constructions. To make the connection more explicit, let us give an equivalent description of weak unitary tensor functors in terms of collections of Hilbert module maps satisfying a system of quadratic relations. 

\smallskip

Assume $G$ is a compact quantum group and $F\colon\Rep G\to\Hilb A$ is a weak unitary tensor functor. As in Section~\ref{s2}, fix representatives $U_\alpha$ of irreducible unitary representations of $G$, and assume $U_e$ is the trivial representation. Consider the correspondences $M_\alpha=F(U_\alpha)$ and the linear maps $\varphi^\gamma_{\alpha,\beta}$ from $\Mor(U_\alpha\times U_\beta,U_\gamma)$ into the space of bounded $A$-bilinear maps $M_\alpha\otimes_AM_\beta\to M_\gamma$ defined by $\varphi^\gamma_{\alpha,\beta}(T)=F(T)F_2$. We then have the following:
\begin{itemize}
\item[(i)] $M_e=A$;
\item[(ii)] if a morphism of the form $(T_1,\dots,T_n)\colon U_\alpha\times U_\beta\to\oplus^n_{i=1}U_{\gamma_i}$ is unitary, then the map $(\varphi(T_1),\dots,\varphi(T_n))\colon M_\alpha\otimes_A M_\beta\to\oplus^n_{i=1}M_{\gamma_i}$ is isometric;
\item[(iii)] the image of the identity map $U_\beta\to U_\beta$ under $\varphi^\beta_{e,\beta}$ is the map $A\otimes_A M_\beta\to M_\beta$ such that $a\otimes X\mapsto aX$, and similarly the image of the identity map $U_\alpha\to U_\alpha$ under $\varphi^\alpha_{\alpha,e}$ is the map $M_\alpha\otimes_AA\to M_\alpha$ such that $X\otimes a\mapsto Xa$;
\item[(iv)] if a morphism $U_\alpha\times U_\beta\times U_\gamma\to U_\delta$ is written as $\sum_iS_i(T_i\otimes\iota)=\sum_jS_j'(\iota\otimes T'_j)$ for some morphisms $T_i\colon U_\alpha\times U_\beta\to U_{\alpha_i}$, $S_i\colon U_{\alpha_i}\times U_\gamma\to U_\delta$, $T'_j\colon U_\beta\times U_\gamma\to U_{\beta_j}$ and $S'_j\colon U_{\alpha}\times U_{\beta_j}\to U_\delta$, then
$$
\sum_i\varphi(S_i)(\varphi(T_i)\otimes\iota)=\sum_j\varphi(S_j')(\iota\otimes\varphi(T'_j)) \ \ \text{as maps}\ \ M_\alpha\otimes_A M_\beta\otimes_AM_\gamma\to M_\eta;
$$
\item[(v)] for every vector $X\in M_\alpha$ and every morphism $T\colon U_\alpha\times U_\beta\to U_\gamma$, the right $A$-linear map $S_X[T]\colon M_\beta\to M_\gamma$ mapping $Y$ into $\varphi(T)(X\otimes Y)$ is adjointable, and if a morphism $U_\gamma\times U_\delta\to U_\alpha\times U_\eta$ is written as $\sum_i(\iota\otimes S_i)(T^*_i\otimes \iota)=\sum_jP_j^*R_j$ for some morphisms $T_i\colon U_{\alpha}\times U_{\beta_i}\to U_{\gamma}$, $S_i\colon U_{\beta_i}\times U_\delta\to U_\eta$, $R_j\colon U_\gamma\times U_\delta\to U_{\gamma_j}$ and $P_j\colon U_{\alpha}\times U_{\eta}\to U_{\gamma_j}$, then
$$
\sum_i\varphi(S_i)(S_X[T_i]^*\otimes\iota)=\sum_jS_X[P_j]^*\varphi(R_j) \ \ \text{as maps}\ \ M_\gamma\otimes_A M_\delta\to M_\eta.
$$
\end{itemize}
Properties (i)-(iv) follow immediately by definition. As will become clear from the proof of the following proposition, the last property, in the presence of the other four, is equivalent to condition~(v) in Definition~\ref{dwutf}. In particular, if in (ii) we have unitary maps instead of isometric maps, then (v) is a consequence of properties (i)-(iv).

\begin{proposition}
Assume we are given correspondences $M_\alpha\in\Hilb A$ and linear maps $\varphi^\gamma_{\alpha,\beta}$ from $\Mor(U_\alpha\times U_\beta,U_\gamma)$ into the space of bounded $A$-bilinear maps $M_\alpha\otimes_AM_\beta\to M_\gamma$ such that the above conditions (i)-(v) are satisfied. Then there exists a unique, up to a natural unitary  monoidal  isomorphism, weak unitary tensor functor $F\colon\Rep G\to\Hilb A$ such that $F(U_\alpha)=M_\alpha$ and $F(T)F_2=\varphi^\gamma_{\alpha,\beta}(T)$ for $T\in\Mor(U_\alpha\times U_\beta,U_\gamma)$.
\end{proposition}

\bp By virtue of semisimplicity of $\Rep G$, there exists a unique, up to a natural unitary isomorphism, unitary functor $F\colon\Rep G\to\Hilb A$ such that $F(U_\alpha)=M_\alpha$. We then define
$$
F_2\colon F(U_\alpha)\otimes_AF(U_\beta)\to F(U_\alpha\times U_\beta) \ \ \text{by}\ \ F_2=\sum_iF(w_i^*)\varphi(w_i),
$$
where $w_i\colon U_\alpha\times U_\beta\to U_{\gamma_i}$ are coisometric morphisms such that $\sum_iw_i^*w_i=\iota$. It is easy to see that this definition does not depend on the choice of $w_i$. By condition (ii) the map $F_2$ is isometric. Note also that $F(w_i)F_2=\varphi(w_i)$, whence $F(T)F_2=\varphi(T)$ for all $T\in\Mor(U_\alpha\times U_\beta,U_\gamma)$. By semisimplicity of $\Rep G$ the isometries $F_2\colon F(U_\alpha)\otimes_AF(U_\beta)\to F(U_\alpha\times U_\beta)$ uniquely define a family of natural isometries $F_2\colon F(U)\otimes_AF(V)\to F(U\times V)$. The only not entirely obvious property left to check is commutativity of two diagrams in Definition~\ref{dwutf}.

\smallskip

For the first diagram, we have to check that $F_2(F_2\otimes\iota)=F_2(\iota\otimes F_2)$ as maps $F(U_\alpha)\otimes_A F(U_\beta)\otimes_AF(U_\gamma)\to F(U_\alpha\times U_\beta\times U_\gamma)$. It suffices to check that $F(w)F_2(F_2\otimes\iota)=F(w)F_2(\iota\otimes F_2)$ for any morphism $w\colon U_\alpha\times U_\beta\times U_\gamma\to U_\delta$. Any such morphism can be written as $\sum_iS_i(T_i\otimes\iota)=\sum_jS_j'(\iota\otimes T'_j)$. By naturality of $F_2$ we then have
$$
F(w)F_2(F_2\otimes\iota)=\sum_iF(S_i)F_2(F(T_i)F_2\otimes\iota) =\sum_i\varphi(S_i)(\varphi(T_i)\otimes\iota),
$$
and similarly $F(w)F_2(\iota\otimes F_2)=\sum_j\varphi(S_j')(\iota\otimes\varphi(T'_j))$. By condition (iv) these expressions are equal.

\smallskip

It remains to show that for every $X\in F(U)$ the maps $S_X=S_{X,V}\colon F(V)\to F(U\otimes V)$ are adjointable and $F_2(S^*_X\otimes\iota)=S_X^*F_2$. For the adjointability it suffices to show that the map $S_X\colon F(U_\beta)\to F(U_\alpha\times U_\beta)$ is adjointable for every $X\in F(U_\alpha)$. Decomposing $U_\alpha\times U_\beta$ into irreducible representations, we see that adjointability of $S_X$ is equivalent to adjointability of $F(T)S_X$ for all morphisms $T\colon U_\alpha\times U_\beta\to U_\gamma$. Since $F(T)S_X(Y)=F(T)F_2(X\otimes Y)=\varphi(T)(X\otimes Y)$, we have $F(T)S_X=S_X[T]$, so adjointability of $F(T)S_X$ is part of condition (v).

Finally, we have to show that $F_2(S^*_X\otimes\iota)=S_X^*F_2$ as maps $F(U_\alpha\times U_\beta)\otimes_A F(U_\delta)\to F(U_\beta\times U_\delta)$ for $X\in F(U_\alpha)$. This is equivalent to $$F(S)F_2(S^*_X\otimes\iota)(F(T^*)\otimes\iota)=F(S)S_X^*F_2(F(T^*)\otimes\iota)\ \ \text{as maps}\ \ F(U_\gamma)\otimes_A F(U_\delta)\to F(U_\eta)$$ for all morphisms $S\colon U_\beta\times U_\delta\to U_\eta$ and $T\colon U_\alpha\times U_\beta\to U_\gamma$. The left hand side of the above identity equals $\varphi(T)(S_X[T]^*\otimes\iota)$, while the right hand side, by \eqref{enats}, equals
$$
S^*_XF(\iota\otimes S)F_2(F(T^*)\otimes\iota)=S^*_XF((\iota\otimes S)(T^*\otimes\iota))F_2.
$$
Writing the morphism $(\iota\otimes S)(T^*\otimes\iota)$ as $\sum_jP_j^*R_j$ for some $R_j\colon U_\gamma\times U_\delta\to U_{\gamma_j}$ and $P_j\colon U_{\alpha}\times U_{\eta}\to U_{\gamma_j}$, we can write the last expression as
$$
\sum_j S^*_XF(P_j^*)F(R_j)F_2=\sum_jS_X[P_j]^*\varphi(R_j).
$$
We thus see that the identity $F_2(S^*_X\otimes\iota)=S_X^*F_2$ follows from condition (v). Furthermore, from the proof we see that it is equivalent to that condition, since any morphism $U_\gamma\times U_\delta\to U_\alpha\times U_\eta$ can be written as $\sum_i(\iota\otimes S_i)(T^*_i\otimes \iota)$ for appropriate morphisms $T_i\colon U_{\alpha}\times U_{\beta_i}\to U_{\gamma}$ and $S_i\colon U_{\beta_i}\times U_\delta\to U_\eta$ using Frobenius reciprocity.
\ep

Given data $\{M_\alpha\}_\alpha$ and $\{\varphi^\gamma_{\alpha,\beta}\}_{\alpha,\beta,\gamma}$ as above, the construction of the corresponding C$^*$-algebra~$B_F$ from Section~\ref{s2} goes as follows. For every $\alpha$ define a new scalar product on $\bar H_\alpha=\bar H_{U_\alpha}$~by
$$
(\bar\xi,\bar\zeta)=\frac{1}{\dim_qU_\alpha}(\zeta,\rho^{-1}\xi).
$$
Consider the right Hilbert $A$-module
$$
M=\ell^2\text{-}\bigoplus_\alpha \bar H_\alpha\otimes M_\alpha.
$$
For every $x=\bar\xi\otimes X\in\bar H_\alpha\otimes M_\alpha$ define an operator $L_x$ on $M$ by
$$
L_x(\bar\zeta\otimes Y)=\sum_i\overline{w_i(\xi\otimes\zeta)}\otimes\varphi^{\gamma_i}_{\alpha,\beta}(w_i)(X\otimes Y)\ \ \text{for}\ \ \bar\zeta\otimes Y\in \bar H_\beta\otimes M_\beta,
$$
where $w_i\in\Mor(U_\alpha\times U_\beta,U_{\gamma_i})$ are coisometries such that $\sum_i w_i^*w_i=\iota$. Then the results of Section~\ref{s2} can be summarized by saying that the operators $L_{\bar\xi\otimes X}$ for all $\xi\in H_\alpha$, $X\in M_\alpha$ and all indices $\alpha$, span a $*$-algebra of bounded operators on the Hilbert $A$-module $M$, and~$B_F$ is the norm closure of this algebra.

\begin{example}
Assume $G$ is the dual of a discrete group $\Gamma$. We identify the set of isomorphism classes of irreducible representations of $G$ with $\Gamma$. Then, up to equivalence, a weak unitary tensor functor $F\colon\Rep G\to\Hilb A$ is the same as a collection of C$^*$-correspondences $M_\alpha$, $\alpha\in\Gamma$, over $A$, together with $A$-bilinear isometries $\varphi_{\alpha,\beta}\colon M_\alpha\otimes_A M_\beta\to M_{\alpha\beta}$ such that
\begin{itemize}
\item[(a)] $M_e=A$;
\item[(b)] $\varphi_{e,\alpha}\colon A\otimes_A M_\alpha\to M_\alpha$ and $\varphi_{\alpha,e}\colon M_\alpha\otimes_A A\to M_\alpha$ are the maps $a\otimes X\mapsto aX$ and $X\otimes a\mapsto Xa$, respectively;
\item[(c)] $\varphi_{\alpha\beta,\gamma}(\varphi_{\alpha,\beta}\otimes\iota) =\varphi_{\alpha,\beta\gamma}(\iota\otimes\varphi_{\beta,\gamma})$;
\item[(d)] for every vector $X\in M_\alpha$ and $\beta\in\Gamma$, the map $S_X\colon M_\beta\to M_{\alpha\beta}$, $Y\mapsto\varphi_{\alpha,\beta}(X\otimes Y)$, is adjointable, and $\varphi_{\beta,\gamma}(S^*_X\otimes\iota)=S^*_X\varphi_{\alpha\beta,\gamma}$ as maps $M_{\alpha\beta}\otimes_A M_\gamma\to M_{\beta\gamma}$.
\end{itemize}
This is similar to the definition of product systems of C$^*$-correspondences~\cite{Fow}. The difference is that instead of semigroups we consider groups, the maps $\varphi_{\alpha,\beta}$ are not assumed to be unitary, but then the additional assumption (d) is required. We remind again that if the maps $\varphi_{\alpha,\beta}$ are unitary, condition~(d) is not needed.

Since conditions (a)-(d) describe spectral subspaces of an arbitrary coaction of $\Gamma$, our results for $G=\hat\Gamma$ simply mean that these conditions give an equivalent characterization of Fell bundles over~$\Gamma$~\cite{FD2}. Explicitly, the $*$-structure on the bundle $\{M_\alpha\}_{\alpha\in\Gamma}$ is given by the operation $\bullet$ defined in Lemma~\ref{lbullet}, so $X^\bullet=S^*_X(1)\in M_{\alpha^{-1}}$ for $X\in M_\alpha$ if $A$ is unital, and in general~$X^\bullet$ is characterized by $\langle X^\bullet,Y\rangle=\varphi_{\alpha,\alpha^{-1}}(X\otimes Y)$ for $Y\in M_{\alpha^{-1}}$. Clearly, $B_F$ is nothing else than the cross-sectional C$^*$-algebra of this bundle.\ee
\end{example}

\begin{example}
Assume $G=\T$. Let $M$ be a Hilbert $A$-bimodule, meaning that $M$ carries the structures of a right Hilbert $A$-module with inner product $\langle\cdot,\cdot\rangle_R$ and of a left Hilbert $A$-module with inner product $_L\langle\cdot,\cdot\rangle$, and $X\langle Y,Z\rangle_R={}_L\langle X,Y\rangle Z$ for all $X,Y,Z\in M$. Consider the complex conjugate Hilbert $A$-bimodule~$\bar M$, so $a\bar X=\overline{Xa^*}$, $\bar Xa=\overline{a^*X}$, $_L\langle \bar X,\bar Y\rangle=\langle X,Y\rangle_R$ and $\langle \bar X,\bar Y\rangle_R={_L}\langle X,Y\rangle$. Define C$^*$-correspondences $M_n$, $n\in\Z$, over $A$ by $M_0=A$, $M_n=M^{\otimes_A n}$ for $n\ge1$ and $M_n=\bar M^{\otimes_A |n|}$ for $n\le-1$. We have obvious isometries $\varphi_{m,n}\colon M_m\otimes_A M_n\to M_{m+n}$. In order to show that they define a weak unitary tensor functor $F\colon\Rep\T\to\Hilb A$, we have to check conditions (a)-(d) from the previous example. Conditions (a) and (b) are obviously satisfied. A moment's reflection shows that since the maps $\varphi_{m,n}$ are surjective for $m$ and $n$ of the same sign, it suffices to check the other two conditions only for $\alpha=\pm1$. For such $\alpha$ conditions (c) and (d) easily follow from the identity $X\langle Y,Z\rangle_R={}_L\langle X,Y\rangle Z$.

The corresponding C$^*$-algebra $B_F$ is the algebra $A\rtimes_M\Z$, the crossed product of $A$ by the Hilbert $A$-bimodule $M$, defined in~\cite{AAE}, where it was shown directly that $\{M_n\}_{n\in\Z}$ forms a Fell bundle over~$\Z$. Recall that the C$^*$-algebra $A\rtimes_M\Z$ is canonically isomorphic to the Cuntz-Pimsner algebra $\mathcal O_M$.\ee
\end{example}

Let us return to the case of a general compact quantum group $G$. Recall that a unitary $2$-cocycle on the dual discrete quantum group $\hat G$ is a unitary element
$$
\Omega\in W^*(G)\bar\otimes W^*(G)\subset(\C[G]\otimes\C[G])^*
$$
such that $(\Omega\otimes1)(\Dhat\otimes\iota)(\Omega)=(1\otimes\Omega)(\iota\otimes\Dhat)(\Omega)$. Any such cocycle defines a unitary fiber functor $E_\Omega\colon\Rep G\to\Hilbf$ that is identity on objects and morphisms, while the tensor structure $E_\Omega(U)\otimes E_\Omega(V)\to E_\Omega(U\times V)$ is given by the action of $\Omega^*$. By Woronowicz's Tannaka-Krein duality,
this functor defines a new deformed compact quantum group $G_\Omega$. More concretely, we have $W^*(\hat G_\Omega)=W^*(\hat G)$ as von Neumann algebras, while the new coproduct is given by $\Dhat_\Omega(\omega)=\Omega\Dhat(\omega)\Omega^*$. Equivalently, $\C[G_\Omega]=\C[G]$ as coalgebras, while the new product and involution are obtained by dualizing $(W^*(G_\Omega),\Dhat_\Omega)$.

Assume now that we have a continuous left action $\theta$ of the reduced form of $G$ on a C$^*$-algebra~$B$ with fixed point algebra $A$. Consider the corresponding spectral functor $F\colon\Rep G\to\Hilb A$. Since by construction the categories $\Rep G$ and $\Rep G_\Omega$ are equivalent, we can consider $F$ as a weak unitary tensor functor $\Rep G_\Omega\to\Hilb A$. It defines a C$^*$-algebra $B_\Omega$ carrying a continuous left action of the reduced form of $G_\Omega$. These algebras were defined and studied in greater generality in~\cite{DRV} and~\cite{NT} (and in various special cases by many authors earlier; in particular, see \cite{Yam} for the case of the dual of a discrete group). But as we will see in a moment, the categorical picture provides a very simple and concrete approach.

In order to formulate the result we need to introduce more notation. First of all we will need a special element $u\in\U(G)=\C[G]^*$ defined by $u=m(\iota\otimes\hat S)(\Omega)$, where $m\colon (\C[G]\otimes\C[G])^*\to\C[G]^*$ is the product map, which is by definition dual to the coproduct on~$\C[G]$. The element $u$ is invertible, with inverse given by
$$
u^{-1}=m(\hat S\otimes\iota)(\Omega^*)=\hat S(u^*).
$$
The antipode on the dual of $G_\Omega$ is given by $\hat S_\Omega=u\hat S(\cdot)u^{-1}$, and correspondingly the involution $\dagger$ on~$\C[G_\Omega]$ is given by
$$
a^\dagger=\big[(u^{-1})^*\otimes\iota\otimes u^*)\Delta^{(2)}(a)\big]^*
=(u^*\otimes\iota\otimes (u^{-1})^*)\Delta^{(2)}(a^*),
$$
see for instance~\cite[Example~2.3.9]{NTbook}. It is easy to check that the element $u$ can also be characterized by the identities
\begin{equation}\label{eu}
\Omega R_U=(u\otimes\iota)R_U\ \ \text{as maps}\ \ \C\to H_{\bar U}\otimes H_U,
\end{equation}
see~\cite[Example~2.2.23]{NTbook}.

Next, consider the subalgebra $\B\subset B$ of regular elements. Then the map $\B\otimes\U(G)\to\B$, $x\otimes\omega\mapsto x\lhd\omega=(\omega\otimes\iota)\theta(x)$, defines a right $\U(G)$-module structure on $\B$.

\begin{proposition}
With the above notation, the following formulas define a new $*$-algebra $\B_\Omega$ with underlying space $\B$, product $\star$ and involution $\dagger$:
$$
x\star y=m((x\otimes y)\lhd\Omega),\ \ x^\dagger=x^*\lhd u^*,
$$
where $m\colon\B\otimes\B\to\B$ is the original product map. Furthermore, the map $\theta$, considered as a map $\B_\Omega\to\C[G_\Omega]\otimes\B_\Omega$, defines a left algebraic action of $G_\Omega$ on $\B_\Omega$.
\end{proposition}

\bp By multiplying $\Omega$ by a phase factor we may assume that $\Omega$ is counital, that is, $(\hat\eps\otimes\iota)(\Omega)=(\iota\otimes\hat\eps)(\Omega)=1$.

For every finite dimensional  unitary representation $U\in B(H_U)\otimes\C[G]$ denote by $U^\Omega$ the same element $U$ considered as an element of $B(H_U)\otimes\C[G_\Omega]$. Then we have a unitary monoidal equivalence of categories $E^\Omega\colon \Rep G_\Omega\to \Rep G$ such that $E^\Omega(U^\Omega)=U$, $E$ is the identity map on morphisms, and $E^\Omega_2\colon E^\Omega(U^\Omega)\otimes E^\Omega(V^\Omega)\to E^\Omega(U^\Omega\times V^\Omega)$ is given by $\Omega\colon H_U\otimes H_V\to H_U\otimes H_V=H_{U\times V}$.

We identify the algebra $\B$ with the algebra $\B_F$ defined by a weak unitary tensor functor $F$. Now we claim that with the above setup the $*$-algebra $\B_{FE^\Omega}$ corresponding to the weak unitary tensor functor $FE^\Omega$ is exactly $\B_\Omega$, and the map $\theta_{FE^\Omega}$ coincides with $\theta$. Note that counitality of $\Omega$ is needed for condition (iii) in Definition~\ref{dwutf} to be satisfied by the functor $FE^\Omega$.

\smallskip

As linear spaces, we have
$$
\B_{FE^\Omega}=\bigoplus_\alpha\bar H_{U^\Omega_\alpha}\otimes FE^\Omega(U^\Omega_\alpha)=\bigoplus_\alpha\bar H_\alpha\otimes F(U_\alpha)=\B.
$$
Denote by $\star$ the product on $\B_{FE^\Omega}$. Note that if $w\colon H_\gamma\to H_\alpha\otimes H_\beta$ is a morphism $U_\gamma\to U_\alpha\times U_\beta$, then $\Omega w$ is a morphism $U^\Omega_\gamma\to U^\Omega_\alpha\times U^\Omega_\beta$. From this we get that if $x=\bar\xi\otimes X\in\bar H_\alpha\otimes F(U_\alpha)$ and $y=\bar\zeta\otimes Y\in\bar H_\beta\otimes F(U_\beta)$, then
$$
x\star y=\sum_i \overline{w_i^*\Omega^*(\xi\otimes\zeta)}\otimes F(w_i^*)F_2(X\otimes Y),
$$
where $w_i\in\Mor(U_{\gamma_i},U_\alpha\times U_\beta)$ are isometries such that $\sum_iw_iw_i^*=\iota$. Since the right $\U(G)$-module structure on $\B=\B_F$ is given by $(\bar\eta\otimes Z)\lhd\omega=\overline{\omega^*\eta}\otimes Z$, the above identity means exactly that
$$
x\star y=m((x\otimes y)\lhd\Omega).
$$

\smallskip

Denote the involution on $\B_{FE^\Omega}$ by $\dagger$. Take $x=\bar\xi\otimes X\in\bar H_\alpha\otimes F(U_\alpha)$. We may assume that $\bar U_\alpha=U_{\bar\alpha}$ for some index $\bar\alpha$. Then, by definition,
$$
x^\dagger=\overline{\xi^\#}\otimes X^\#,
$$
where $\xi^\#\in H_{\bar\alpha}$ is such that $(R^{\Omega}_\alpha)^*(\zeta\otimes\xi)=(\zeta,\xi^\#)$ for all $\zeta\in H_{\bar\alpha}$ and $X^\#\in F(U_{\bar\alpha})$ is such that $(FE^\Omega)(\bar R^{\Omega}_\alpha)^*(FE^\Omega)_2(X\otimes Y)=\langle X^\#,Y\rangle$ for all $Y\in F(U_{\bar\alpha})$, where $R^\Omega_\alpha$ and $\bar R^\Omega_\alpha$ solve the conjugate equations for $U^\Omega_\alpha$ and $U^\Omega_{\bar\alpha}$. Note that by irreducubility the operation $\#$ depends on the choice of such a solution, but $\overline{\xi^\#}\otimes X^\#$ does not. Taking the solutions $R_\alpha$ and $\bar R_\alpha$ of the conjugate equations for $U_\alpha$ and $U_{\bar\alpha}$ defined by \eqref{econj}, we can take $R^\Omega_\alpha=\Omega R_\alpha$ and $\bar R^\Omega_\alpha=\Omega\bar R_\alpha$. In this case $X^\#=X^\bullet$, while for $\xi^\dagger$, using~\eqref{eu}, we get
$$
(R^{\Omega}_\alpha)^*(\zeta\otimes\xi)=R^*_\alpha(u^*\zeta\otimes\xi)=(u^*\zeta,\xi^\bullet),
$$
so $\xi^\#=u\xi^\bullet$. Therefore
$$
x^\dagger=\overline{u\xi^\bullet}\otimes X^\bullet=(\overline{\xi^\bullet}\otimes X^\bullet)\lhd u^*=x^*\lhd u^*.
$$

\smallskip

Finally, the maps $\theta=\theta_F$ and $\theta_{FE^\Omega}$ coincide on $\B_{FE^\Omega}$, since they both define the same right $\U(G)$-module structure, given by $(\bar\eta\otimes Z)\lhd\omega=\overline{\omega^*\eta}\otimes Z$.
\ep

The construction of a new product on a module algebra using a cocycle on a Hopf algebra is, of course, well-known. The point of the above proposition is that it effortlessly gives not only the new product, but also the $*$-structure and the existence of a unique C$^*$-completion of the algebra carrying an action of the reduced deformed quantum group.

\bigskip

\end{document}